\documentclass[12pt]{amsart}
\usepackage{amssymb}
\usepackage{epsfig}
\usepackage{graphicx}
\numberwithin{equation}{section}

\input xy
\xyoption{all}

\advance\headheight 2pt
\calclayout
\allowdisplaybreaks[3]

\theoremstyle{plain}
\newtheorem{prop}{Proposition}

\newtheorem{theo}[prop]{Theorem}
\newtheorem{coro}[prop]{Corollary}

\newtheorem{lemm}[prop]{Lemma}

\theoremstyle{definition}
\newtheorem{defi}[prop]{Definition}

\newtheorem{conj}[prop]{Conjecture}

\newtheorem{rema}[prop]{Remark}

\newtheorem{exam}[prop]{Example}

\newcommand{\bA}{\mathbb A}

\newcommand{\bP}{\mathbb P}
\newcommand{\bQ}{\mathbb Q}
\newcommand{\bZ}{\mathbb Z}

\newcommand{\cN}{\mathcal N}
\newcommand{\cO}{\mathcal O}
\newcommand{\cT}{\mathcal T}
\newcommand{\cU}{\mathcal U}

\newcommand{\cX}{\mathcal X}
\newcommand{\cY}{\mathcal Y}

\newcommand{\ra}{\rightarrow}

\newcommand{\Spec}{{\rm Spec}}

\newcommand{\Hilb}{{\mathcal H}\!{\mathit ilb}}

\newcommand{\mfm}{\mathfrak m}
\newcommand{\mfq}{\mathfrak q}

\newcommand{\mfw}{\mathfrak w}
\newcommand{\mfo}{\mathfrak o}

\makeatother
\makeatletter

\author{Brendan Hassett}
\address{Department of Mathematics \\
Rice University, MS 136 \\
Houston, Texas  77251-1892 \\
USA}
\email{hassett@rice.edu}

\author{Yuri Tschinkel}
\address{
                Mathematisches Institut \\
                Bunsenstr. 3-5  \\
        37073 G\"ottingen \\
                Germany \\
Courant Institute, NYU \\ 251 Mercer str. \\ New York, NY 10012 \\ USA}
\email{
tschinkel@cims.nyu.edu}

\title[Approximation at places of bad reduction]{Approximation 
at places of bad reduction for rationally connected varieties}

\begin{document}
\date{\today}

\begin{abstract}
This paper addresses weak approximation for rationally connected
varieties defined over the function field of a curve, especially
at places of bad reduction.  Our
approach entails analyzing the rational connectivity of the smooth 
locus of singular reductions of the variety.  As an application, 
we prove weak approximation for cubic surfaces with square-free
discriminant.  
\end{abstract}

\maketitle
\tableofcontents

\section{Introduction}
\label{sect:introduction}

In number theory, many results and techniques rely on 
approximating adelic points by rational points. 
In this paper, we study geometric versions of these notions 
for rationally connected varieties over the function field of a curve.
In this context, rational points correspond to sections of 
rationally-connected fibrations over the curve.  
We are looking for sections with prescribed jet data 
in finitely many fibers 
(see Section~\ref{sect:approx} for definitions). 

Let $k$ be an algebraically closed field of characteristic zero,
$B$ a smooth curve over $k$ with function field $F=k(B)$. 
Let $\overline{B}$ be the smooth projective model of $F$
and put $S:=\overline{B}\setminus B$.

\begin{theo}\label{theo:main}
Let $X$ be a smooth proper rationally connected variety over $F$, and $\pi:\cX \ra B$ a
model of $X$, i.e., $\cX$ is an algebraic space flat and proper over $B$
with generic fiber $X$.  Let $\cX^{sm}$ be the locus where $\pi$ is smooth and
$\cX^{\bullet} \subset \cX^{sm}$ be such that
\begin{enumerate}
\item{
there exists a section $s:B\ra \cX^{\bullet}$;}
\item{
for each $b\in B$ and $x\in \cX^{\bullet}_b$, there exists a rational
curve $f:\bP^1\ra \cX^{\bullet}_b$ containing $x$ and the generic point of $\cX^{\bullet}_b$.}
\end{enumerate}
Then sections of $\cX^{\bullet}\ra B$ satisfy approximation away from $S$. 
\end{theo}
Rationally-connected fibrations over curves have sections
by \cite{GHS}.  The existence of a section
through a finite set of prescribed points is addressed
\cite{KMM} 2.13 and \cite{kollar} IV.6.10.1.
Weak approximation is known in fibers of good reduction \cite{HT}, so we 
take simultaneous resolutions of singular fibers of $\cX$ whenever possible
\cite{Br1} \cite{Br2}.  Consequently, when $\cX \ra B$ admits a simultanteous resolution over
some \'etale neighborhood of $b$, we replace $\cX$ by this resolution.
However, the resolved family may be an algebraic space, rather than
a scheme, over $B$.  
This is why Theorem~\ref{theo:main} is stated in this generality.  

We shall actually prove a stronger result, Theorem~\ref{theo:mainbis},
which is applicable in positive characteristic.  In this context,
Corollary~\ref{coro:positive} gives weak approximation at places of 
good reduction.

\

There are very few instances where weak approximation over function fields is known
at all places \cite{CT}:
\begin{itemize}
\item stably rational varieties;
\item connected linear algebraic groups and homogeneous spaces
for these groups;
\item homogeneous space fibrations over varieties that satisfy weak approximation,
for example, conic bundles over rational varieties;
\item Del Pezzo surfaces of degree at least four.
\end{itemize}
Even the case of cubic surfaces remains open, in general.
Madore established weak approximation for cubic surfaces
at places of good reduction \cite{Ma}.  His proof uses the abundance of
distinct unirational parametrizations, and builds on ideas of Swinnerton-Dyer \cite{SWD}.  

When is Theorem~\ref{theo:main} applicable?  
Let $X$ be a smooth projective rationally connected variety over $F=k(B)$,
with $B$ projective.  There
exists a regular model $\pi:\cX \ra B$, and any section $s:B\ra \cX$
is contained in $\cX^{sm}$.   
For each singular fiber $\cX_b$, fix an irreducible component
$\cX^{\bullet}_b\subset \cX_b^{sm}$;  these determine
an open subset $\cX^{\bullet}\subset \cX^{sm}$.  To prove weak
approximation for $X$, it suffices to prove approximation
for each $\cX^{\bullet}$ obtained in this way.  
We do not know how to verify (1) in general:  Is there 
{\em any} section meeting a prescribed irreducible component of $\cX^{sm}_b$?
Further, there is no general result giving a regular model $\cX\ra B$ such that
each irreducible component of $\cX^{sm}_b$ has the property (2).

We give applications to cubic surfaces:
\begin{theo}\label{theo:cubic}
Let $X$ be a smooth cubic surface over $F$ and $\pi:\cX \ra B$ a model
whose singular fibers are cubic surfaces with rational double points.  Suppose
there exists a section $s:B\ra \cX^{sm}$.
Then sections of $\cX^{sm}\ra B$ satisfy approximation away from $S$.  
\end{theo}
When the model is regular all sections are contained in the smooth locus,
so we conclude:
\begin{coro}\label{coro:cubicsm}
Let $X$ be a smooth cubic surface over $F$.  Suppose $X$ admits a regular model $\pi:\cX \ra B$
whose singular fibers are cubic surfaces with rational double points.  
Then weak approximation holds for $X$ away from $S$.  
\end{coro}

There exist cubic surfaces which do not admit models with at most
rational double points in a given fiber, e.g.,
the isotrivial family
$$x^3+y^3+z^3=tw^3$$
over the $t$-line.
Nonetheless, Corollary~\ref{coro:cubicsm} proves weak approximation
for `generic' cubic surfaces.  
\begin{coro}\label{coro:cubicgen}
Let $\Hilb=\bP(\Gamma(\cO_{\bP^3}(3)))\simeq \bP^{19}$ denote the Hilbert scheme of 
cubic surfaces, $\cU\ra \Hilb$ the universal family, 
$D\subset \Hilb$ the discriminant divisor, and $B\subset\Hilb$
a smooth curve transverse to $D$ (i.e., the discriminant is square-free along $B$).  
Then sections of 
$$
\cX=\cU\times_{\Hilb}\Spec(F)\ra B
$$
satisfy approximation away from $S$.  
\end{coro}

Meeting the discriminant transversally is an open condition
on the classifying map to the Hilbert scheme.  
The transversality implies that near singular points 
of $\cX_b$, the model $\cX:=\cU\times_{\Hilb}B$ has local analytic equation
$x^2+y^2+z^2=t$, where $t$ is a local uniformizer for $B$ at $b$.  
In particular, $\cX$ is a regular model and Corollary~\ref{coro:cubicsm}
applies.

\

In our approach to approximation, we require precise control over proper rational curves in
the smooth locus.  
One focus of this paper is to extend standard results on smooth proper rationally
connected varieties to the non-proper case (see Section~\ref{sect:SRC}). 
The application to cubic surfaces entails a refinement of rational connectivity results 
of \cite{KeMc} (see Section~\ref{sect:cubic}).

\

\noindent
{\bf Acknowledgments:}  We are grateful to J. L. Colliot-Th\'el\`ene for numerous discussions about
the problems considered here;  the ideas here were developed during visits to
Orsay by both authors.  We also benefitted from conversations with S. Keel,
A. Knecht, J. Koll\'ar, and J. McKernan.  
The first author was partially supported by the
Sloan Foundation and NSF Grants 0134259 and 0196187.

\section{Notions of approximation}
\label{sect:approx}
Let $F$ be a global field, i.e., a number field or the function
field of a curve $B$ defined over an algebraically closed field $k$.
Let $S$ a finite
set of places of $F$ containing the archimedean places, $\mfo_{F,S}$
the corresponding ring of integers, and   
$\bA_{F,S}$ the restricted direct product over all places outside $S$. 

Let $X$ be an algebraic variety over $F$, 
$X(F)$ the set of $F$-rational points
and $X(\bA_{F,S})\subset \prod_{v\notin S}X(F_v)$ the set of $\bA_{F,S}$-points of $X$.
The set $X(\bA_{F,S})$ 
carries a natural direct product topology.
One says that {\em weak approximation} holds for $X$ away from $S$
if $X(F)$ is dense in this topology.  

The set $X(\bA_{F,S})$ also carries a natural adelic topology:
The basic open subsets are 
$$\prod_{v\in S'}\mathfrak u_v\times \prod_{v\notin (S\cup S')} \cX(\mfo_v),$$
where $S'$ is a finite set of nonarchimedean places disjoint from $S$,
$\cX$ is a model over $\Spec(\mfo_{F,S})$,
$\mfo_v$ is the completion of $\mfo_{F,S}$ at $v$, and 
$\mathfrak u_v \subset X(F_v)$ an open subset in the $v$-adic
analytic topology on $X(F_v)$.  This does not depend on the choice of model.
{\em Strong approximation} holds for $X$ away from $S$
if $X(F)$ is dense in $X(\bA_{F,S})$. 
Note that strong approximation 
implies weak approximation. 
Conversely, for
$\cX$ proper over the integers, weak approximation implies
strong approximation, since $\cX(\mfo_v)=X(F_v)$;  in these
cases, we will use the term
weak approximation for the sake of consistency.

Finally, there is a formulation which
is sensitive to the choice of model $\cX$.  
Consider the topology on $\prod_{v\notin S}\cX(\mfo_v)$ with basic open subsets
$$\prod_{v\in S'}\mathfrak u_v\times \prod_{v\notin (S\cup S')} \cX(\mfo_v),$$
with $\mathfrak u_v \subset \cX(\mfo_v)$ an open subset.  
We say that {\em approximation holds for $S$-integral points of $\cX$}
if $\cX(\mfo_{F,S})$ is dense in this product.  This
is a weak version of strong approximation.  

\

We now focus on the function field case:  Let $\overline{B}$ be a smooth 
projective model of $B$ with $S=\overline{B}\setminus B$;
place $v$ correspond to points $b\in \overline{B}$.  
Let $X$ be a smooth variety proper over $F=k(B)$, 
$\pi:\cX \ra B$ a model proper
and flat over $B$ (which exists by \cite{Na}),
and $\cX^{\bullet}\subset \cX^{sm}$ a
model for $X$ surjecting onto $B$.  
Since $\pi$ is proper, $F$-rational points of $X$ correspond to 
sections $s:B\ra \cX$.  If $\cX$ is regular $s$ factors through
$\cX^{sm}$.  

\begin{defi}
An {\em admissible section} of $\pi:\cX \ra B$ is a section
$s:B\ra \cX^{sm}$.
An {\em admissible $N$-jet} of
$\pi$ at $b$ is a section of
$$
\cX^{sm}\times_B \Spec(\cO_{B,b}/\mfm^{N+1}_{B,b})\ra
\Spec(\cO_{B,b}/\mfm^{N+1}_{B,b}).
$$
An {\em approximable $N$-jet} of
$\pi$ at $b$ is a section of
$$
\cX\times_B \Spec(\cO_{B,b}/\mfm^{N+1}_{B,b})\ra
\Spec(\cO_{B,b}/\mfm^{N+1}_{B,b})
$$
that may be lifted to a section of 
$\widehat{\cX}_b \ra \widehat{B}_b,$
with $\widehat{B}_b=\Spec(\hat{\cO}_{B,b})$ and $\widehat{\cX}_b=X\times_B \widehat{B}_b.$
\end{defi}
Hensel's lemma guarantees that every admissible $N$-jet is
approximable.  Let $\{b_i\}_{i\in I}$ be a finite set of points 
and $j_i$ an admissible $N$-jet of $\pi$ at $b_i$.
We write $J=\{ j_i\}_{i\in I}$ for the
corresponding collection of admissible $N$-jets.

The notions of weak and strong approximation
introduced above have geometric interpretations
\begin{itemize}
\item{Weak and strong approximation hold for $X$ away from $S$ if
any finite collection of approximable jets of $\pi$ can be realized
by a section $s:B\ra \cX$.}
\item{This is equivalent to weak 
approximation holding for $X^{\bullet}$ away from $S$:  Every
jet in $\cX$ at $b$ can be realized by a section
$\cX\times_B \widehat{B}_b \ra \widehat{B}_b$ meeting 
$\widehat{\cX}^{\bullet}_b$.}
\item{If $\cX$ is regular these are equivalent to the condition that
any collection of admissible jets of $\pi$ can be realized by
a section $s:B\ra \cX^{sm}$.}
\end{itemize}
There is an analogous formulation of approximation for integral points:
\begin{itemize}
\item{Approximation holds for sections of $\cX^{\bullet}\ra B$ away from $S$
if each collection of  
jet data in $\cX^{\bullet}$ can be realized by a section $s:B\ra \cX^{\bullet}$.}
\item{If $\cX$ is regular and $\cX^{\bullet}=\cX^{sm}$ this is equivalent to
weak approximation for $X$.}
\end{itemize}

\section{Curves, combs, and deformations}
\label{sect:term}

The {\em dual graph} associated with a nodal curve $C$ has
vertices are indexed by the irreducible components of $C$
and its edges indexed by
the intersections of these components.
A projective nodal curve $C$ is {\em tree-like} if
\begin{itemize}
\item{each irreducible component of $C$ is smooth;}
\item{the dual graph of $C$ is a tree.}
\end{itemize}

\begin{defi}
A {\em comb with $m$ reducible teeth} is a projective nodal curve $C$ with
$m+1$ subcurves
$D, T_1, \ldots, T_m$ such that
\begin{itemize}
\item{$D$ is smooth and irreducible;}
\item{$T_l\cap T_{l'}=\emptyset$, for all $l\neq l'$;}
\item{each $T_l$
meets $D$ transversally
in a single point; and }
\item{each $T_l$ is a chain of $\bP^1$'s.}
\end{itemize}
Here $D$ is called the {\em handle} and the
$T_l$ the {\em reducible teeth}.
\end{defi}

We will use the following lemma, which has the same proof as 
Proposition 24 of \cite{HT}:
\begin{lemm}\label{lemm:old}
Let $C$ be a tree-like curve, $W$ a smooth algebraic space,
$h:C\ra W$ an immersion with nodal image.  Suppose that for each irreducible
component $C_l$ of $C$, $H^1(C_l,\cN_h\otimes \cO_{C_l})=0$ and
$\cN_h \otimes \cO_{C_l}$ is globally generated.  Then $h$ deforms to
an immersion.  

Suppose furthermore that $\mfw=\{w_1,\ldots,w_M\} \subset C$ is a collection of 
smooth points such that for each component $C_l$, 
$H^1(\cN_h\otimes \cO_{C_l}(-\mfw))=0$ and the sheaf $\cN_h \otimes \cO_{C_l}(-\mfw)$
admits a section nonzero at each point of the quotient
$$(\cN_h\otimes \cO_{C_l})/\cN_{h|C_l}.$$
Then $h:C\ra W$ deforms to an immersion of a smooth curve into $W$ containing 
$h(\mfw)$.  
\end{lemm}

\section{Strong rational connectivity}
\label{sect:SRC}
\begin{defi}
A variety $X$ is {\em rationally connected} (resp. {\em 
separably rationally connected}) if there is a family 
of proper irreducible rational curves $g:U\ra Z$ (resp.
$\pi_2:U=\bP^1\times Z\ra Z$) and a cycle
morphism $u:U\ra X$ such that
$$
u^{2}:U\times_Z U \ra X\times X$$
is dominant (resp. smooth over the generic point)).
\end{defi}
Intuitively, two generic points of $X$ can be joined by
an irreducible projective rational curve.  Over fields
of characteristic zero, rational connected varieties
are also separably rationally connected \cite{kollar} IV.3.3.1.

The notion of rational connectedness is a bit subtle over countable fields:
For convenience, we work over an uncountable algebraically closed
field.  
Over such a field, rational connectivity is equivalent to
the condition that two very general points of $X$ can
be joined by such a rational curve.  

\begin{defi}
Let $X$ be a smooth algebraic space of dimension $d$ and
$f:\bP^1 \ra X$ a nonconstant morphism, so we have an isomorphism
$$
f^*\cT_X\simeq \cO_{\bP^1}(a_1) \oplus \ldots \oplus \cO_{\bP^1}(a_d)
$$
for suitable integers $a_1,\ldots,a_d$.
Then $f$ is {\em free} (resp. {\em very free}) if each $a_i\ge 0$
(resp. $a_i\ge 1$).
\end{defi}
We refer the reader to \cite{kollar} IV.3 for further facts
about rationally connected varieties.

One technical result will play a prominent r\^ole in our analysis.
\begin{prop}[\cite{kollar} IV.3.9.4] \label{prop:tech}
Let $V$ be a smooth separably rationally connected (not necessarily
proper) variety.  Then there
exists a nonempty subset $V^0\subset V$ characterized
as the largest open subset such that if $v_1,\ldots,v_m\in V^0$
are distinct closed points, then there is a very free curve
in $V^0$ containing these as smooth points.   Moreover, any
rational curve $C\subset V$ that meets $V^0$ is contained
in $V^0$.
\end{prop}
No example where $V^0\neq V$ is known.  
\begin{rema} \label{rema:include}
Let $V_2$ be a smooth variety, $V_1\subset V_2$ a rationally
connected dense open subvariety, and $V_2^0\subset V_2$
the largest open set satisfying the conditions
of Proposition~\ref{prop:tech} .  Then $V_1^0 \subset V_2^0$.  
Thus a point $v\in V_2$ is in $V_2^0$ provided there is 
a rational curve $f:\bP^1\ra V_2$ through $v$ and meeting $V_1^0$.
\end{rema}

\begin{prop} \label{prop:blowup}
Let $V$ be a smooth separably rationally connected variety, and
$\beta:W\ra V$ an iterated blow-up of $V$ along smooth subvarieties.  
Then $\beta^{-1}(V^0)=W^0$.  
\end{prop}
\begin{proof}  
The inclusion $W^0 \subset \beta^{-1}(V^0)$ is straightforward:
Given points $w_1,\ldots,w_m \in W^0$, there is a very free curve
$g:\bP^1 \ra W^0$ containing them;  we may choose this
to be transversal to the exceptional divisor of $\beta$.  
The inclusion of sheaves
$$\cT_W \hookrightarrow \beta^*\cT_V$$
remains an inclusion after pull-back via $g$, as the support 
of the cokernel does not contain $g(\bP^1)$.  The positivity
of $g^*\cT_W$ implies the positivity of $(\beta\circ g)^*\cT_V$,
which means that $\beta\circ g:\bP^1 \ra V$ is also very free.  

For the reverse direction, we may restrict to the case where $W$
is the blow-up of $V$ along a smooth subvariety $Z$ of
codimension $r>1$, with exceptional divisor $E$.  It is
clear that $\beta^{-1}(V^0\setminus Z) \subset W_0$, so 
consider some $w\in \beta^{-1}(z)$ with $z\in Z\cap V^0$.  
It suffices to construct a rational curve containing $w$ and
the generic point of $W$.  

There exists a very free curve $f':\bP^1\ra V^0$ with the
following properties:
\begin{enumerate}
\item{$f'(\bP^1)$ meets $Z$ only at $z$ (we can always deform
a very free curve so that it misses a codimension $\ge 2$ subset);}
\item{$f'(\bP^1)$ is smooth at $z$ and transverse to $Z$.}
\end{enumerate}
Let $g':\bP^1 \ra W$ denote
the lift to $W$, which is free in $W$,
and $w'=g'(0)$.  If $w'=w$ then we are done. 
Otherwise, let $\ell \subset \beta^{-1}(z)\simeq \bP^{r-1}$
denote the line joining $w$ and $w'$.  Since $g'$ is free,
it admits a small deformation to a free curve $g'':\bP^1 \ra W$
with $w'':=g''(0)\in \ell, w''\neq w'$.  (See Figure~\ref{figa}.)

\begin{figure}[h]
\centerline{\hskip 2cm\includegraphics{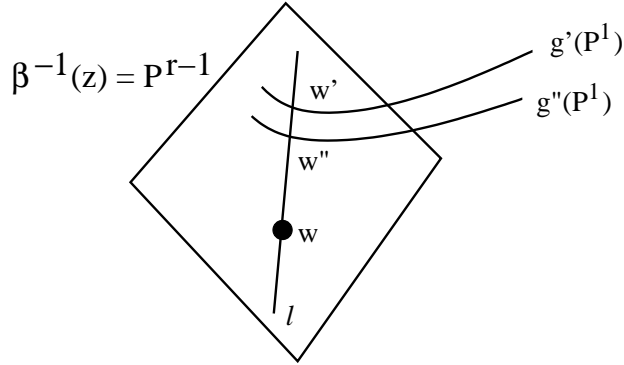}}
\caption{Constructing the comb}
\label{figa}
\end{figure}

We construct a comb
$h:C \ra W$ with handle $\ell\subset \bP^{r-1}\subset W$ and
two teeth $g',g'':\bP^1 \ra W$.  Using the exact sequence
of normal bundles
$$0 \ra \cN_{\ell/E} \ra \cN_{\ell/W} \ra \cN_{E/X}\otimes \cO_{\ell}  \ra 0$$
we find
$$\cN_{\ell/W}\simeq \cO_{\bP^1}^{\dim(V)-r} \oplus 
\cO_{\bP^1}(1)^{r-2}\oplus \cO_{\bP^1}(-1)$$
where the negative summand is in the normal direction to $E$.  
Since $g'(\bP^1)$ and $g''(\bP^1)$ are transverse to $E$, applying
Proposition 23 of \cite{HT} we see
$$\cN_h\otimes \cO_{\ell} \simeq \cO_{\bP^1}^{\dim(V)-r} \oplus 
\cO_{\bP^1}(1)^{r-2}\oplus \cO_{\bP^1}(1);$$
the quotient $(\cN_h \otimes \cO_{\ell}) / \cN_{h|\ell}$
lies in the image of the positive summands.  

Lemma~\ref{lemm:old} implies that
$h:C\ra W$ admits a deformation to a rational curve
containing $w$.  
\end{proof}

A similar argument gives the following 
strengthening of Proposition~\ref{prop:tech} 
(cf. Theorem 2.2 of \cite{De})
\begin{prop} \label{prop:jetstech}
Let $V$ be a smooth separably rationally connected variety and 
$V^0\subset V$ be the distinguished open subset characterized
in Proposition~\ref{prop:tech}.  
Then for any finite collection of jets
$$j_i:\Spec k[\epsilon]/\left<\epsilon^{N+1}\right> \hookrightarrow V^0,
\quad  i=1,\ldots,m$$
supported at distinct points $v_1,\ldots,v_m$,
there exists a very free rational curve 
smooth at $v_1,\ldots,v_m$ with the prescribed jets.
\end{prop}
\begin{proof}
There is an iterated blow-up 
$$\beta:W=W_N \ra \ldots \ra W_j \ra \ldots \ra W_1 \ra V$$ 
and points
$w_1,\ldots,w_m \in W$ so that
if $g:C\ra W$ is a morphism whose image contains $w_1,\ldots,w_m$ then
the image of $f:=\beta\circ g:C\ra V$ contains the given collection
of jets.  
Here is the description:  Over each point $v_i$, we blow up $V$
successively at $N$ points.  Given any smooth curve germ $C$ 
with the prescribed $N$-jet at $v_i$, $W_j$ is the blowup of
$W_{j-1}$ at the points of the proper transform of $C$ lying over the $v_i$.  
Proposition~\ref{prop:blowup} then implies there exists
a very free curve $g:\bP^1 \ra W$ through $w_1,\ldots w_m$.  However,
the image of this curve in $V$ will be singular at $v_i$ if $g(\bP^1)$
meets $\beta^{-1}(v_i)$ in more than one point.  

We claim there exists a very free curve $g_i:\bP^1 \ra W$ meeting 
$\beta^{-1}(v_i)$ only at $w_i$,  transversally.  We choose this
curve so that it is disjoint from $\beta^{-1}(v_j)$ when $j\neq i$.   
Fix generic points $x_i \in g_i(\bP^1)$ and let $g_0:\bP^1 \ra W$ be a very free
curve intersecting $g_i(\bP^1)$ transversely at $x_i$ but not meeting any $\beta^{-1}(v_i)$.
(For example, take $g_0=(\beta^{-1}\circ f_0)$, where $f_0:\bP^1\ra V$ is
a very free curve through $\beta(x_1),\ldots,\beta(x_m)$.)    Consider the
comb $h:C \ra W$ with handle $g_0(\bP^1)$ and $m$-teeth $g_i(\bP^1)$.
This deforms to a very free curve $h':\bP^1 \ra W$ meeting each $\beta^{-1}(v_i)$
only at $w_i$, transversally.  

The proof of the claim is a refinement of the argument for Proposition~\ref{prop:blowup}.
We proceed by induction on $N$.  The base case $N=1$ is contained in
the proof of Proposition ~\ref{prop:blowup}, which gives a very free curve smooth 
at $v_i$ with prescribed tangency.  
Let $E_{i,N}\simeq \bP^{\dim(V)-1}$ be the
last exceptional divisor of $\beta:W \ra V$ over $v_i$, i.e., the exceptional
divisor of the $N$-th blow-up.  For $1\leq j<N$,
let $E_{i,j}\subset W_N$ denote the proper transform of the exceptional 
divisor of $W_j \ra W_{j-1}$ over $v_i$;  we have 
$E_{i,j}\simeq \mathrm{Bl}_{w_{i,j}}\bP^{\dim(V)-1}$, where 
$w_{i,j}$ is the intersection of the proper transform of $C$ with the exceptional
divisor of $W_j \ra W_{j-1}$.  

Suppose that $g'_i:\bP^1 \ra W$ is a very free curve
such that $\beta\circ g'_i$ is smooth with the desired $(N-1)$-jet at $v_i$.  Let
$w'_i=g'_i(\bP^1)\cap \beta^{-1}(v_i)$ denote the unique point of intersection, which we assume is
distinct from $w_i$.  Let $\ell_N$ denote the line in $E_{i,N}\simeq \bP^{\dim(V)-1}$
joining $w_i$ and $w'_i$, and $z_{N-1}$ its point of intersection with $E_{i,N-1}$.  
Let $\ell_{N-1}\subset E_{i,N-1}\simeq \mathrm{Bl}_{w_{i,N-1}}\bP^{\dim(V)-1}$ denote the 
proper transform of a line containing $z_{N-1}$, and $z_{N-2}$ its point of
intersection with $E_{i,N-2}$.  Continue in this way, until we obtain
$\ell_1\subset E_{i,1}$, the proper transform of a line containing $z_1$.  
Finally, let $g''_i:\bP^1 \ra W$ be a very free curve meeting the exceptional
locus transversally at a generic point of $\ell_1$.   (See Figure~\ref{figb}.)
\begin{figure}[h]
\centerline{\hskip 2cm\includegraphics{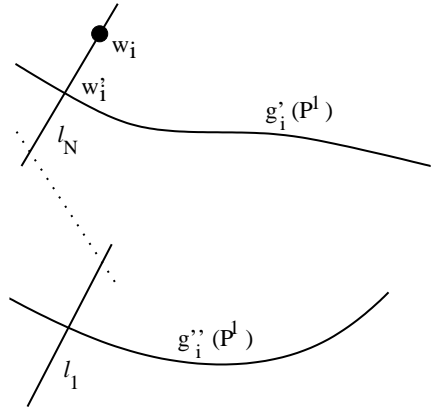}}
\caption{Constructing the comb with reducible teeth}
\label{figb}
\end{figure}

Let $h:C\ra W$ be the comb with handle $\ell_N$ and two reducible teeth:
\begin{enumerate}
\item{$g'_i:\bP^1 \ra W$;}
\item{the union of the lines $\ell_{N-1},\ldots,\ell_1$ and the curve $g''_i:\bP^1 \ra W$;}
\end{enumerate}
By a normal bundle computation similar to that of Proposition~\ref{prop:blowup}, we find that
$\cN_h|\ell_N$ is ample and $\cN_h$ is nonnegative on each of the remaining components:
Again, Lemma~\ref{lemm:old} (or Proposition 24 of \cite{HT}) implies that $h$ admits 
a deformation to an immersed rational curve containg $w_i$.

Here are the details of the computations (cf. \cite{HT} Section 5):
The normal bundle of a line in projective space is
$$\cN_{\ell_N/E_{i,N}}=\cN_{\ell_N/\bP^{\dim(V)-1}}\simeq \cO_{\bP^1}(+1)^{\dim(V)-2}$$
and the normal bundle for an exceptional divisor is
$$\cN_{E_{i,N}/W}\simeq \cO_{\bP^{\dim(V)-1}}(-1).$$
For each $j$ we have
\begin{equation} \label{ES1}
0 \ra \cN_{\ell_j/E_{i,j}} \ra \cN_{\ell_j/W} \ra \cN_{E_{i,j}/W}|_{\ell_j} \ra 0
\end{equation}
which for $j=N$ yields
$$\cN_{\ell_N/W}\simeq \cO_{\bP^1}(+1)^{\dim(V)-2} \oplus \cO_{\bP^1}(-1),$$
with the negative component in the direction normal to $E_{i,N}$.  We also have
an extension
\begin{equation} \label{ES2}
0 \ra \cN_{\ell_j/W} \ra \cN_h|_{\ell_j} \ra Q(\ell_j) \ra 0,
\end{equation} 
where $Q(\ell_j)$ is a torsion sheaf supported at the points where
$\ell_j$ meets the adjacent components.  For $j=N$ these 
are $g_i'(\bP^1)$ and $\ell_{N-1}$, and since the tangent
vectors to these curves are normal to $E_{i,N}$, we find
$$\cN_h|_{\ell_N}\simeq \cO_{\bP^1}(+1)^{\dim(V)-2} \oplus \cO_{\bP^1}(+1).$$

The normal bundle of the proper transform of a line in the blow-up
of projective space at a point of the line is
$$\cN_{\ell_j/E_{i,j}}=\cN_{\ell_j/\mathrm{Bl}_{w_{i,j}}\bP^{\dim(V)-1}}\simeq 
\cO_{\bP^1}^{\dim(V)-2}$$
for $j=1,\ldots,N-1$.  Similarly, we can compute
$$\cN_{E_{i,h}/W}|_{\ell_j}=\cO_{\bP^1}(-2)$$
so the exact sequence analogous to (\ref{ES1}) yields
$$\cN_{\ell_j/W}\simeq \cO_{\bP^1}^{\dim(V)-2} \oplus \cO_{\bP^1}(-2),$$
with the negative component in the direction normal to $E_{i,j}$.  Using
(\ref{ES2}) and the fact that $\ell_j$ is adjacent to $\ell_{j+1}$
and $\ell_{j-1}$ (or $g_i''(\bP^1)$ when $j=1$), we find
$$\cN_h|_{\ell_j}\simeq \cO_{\bP^1}^{\dim(V)-2} \oplus \cO_{\bP^1}.$$
\end{proof}

\begin{defi}
A smooth separably rationally connected 
variety $Y$ is {\em strongly rationally connected} if any
of the following conditions hold:
\begin{enumerate}
\item{for each point $y\in Y$, there exists a rational curve $f:\bP^1 \ra Y$
joining $y$ and a generic point in $Y$;}
\item{for each point $y\in Y$, there exists a free rational curve 
containing $y$;}
\item{for any finite collection of points $y_1,\ldots,y_m\in Y$, there
exists a very free rational curve containing the $y_j$ as smooth points;}
\item{for any finite collection of jets
$$\Spec\, k[\epsilon]/\left<\epsilon^{N+1}\right> \subset Y, i=1,\ldots,m$$
supported at distinct points $y_1,\ldots,y_m$,
there exists a very free rational curve 
smooth at $y_1,\ldots,y_m$ and containing the prescribed jets.}
\end{enumerate}
\end{defi}
The implications
$$(4) \Rightarrow (3) \Rightarrow (2) \Rightarrow (1) $$
are obvious.  
By Proposition~\ref{prop:tech}, 
assertions (1)-(3) are each equivalent to the condition 
$Y=Y^0$.  
Property (4) is analogous to Theorem 2.2 of \cite{De}, which is stated
for proper varieties.  It follows from (1) by 
Proposition~\ref{prop:jetstech}.  

\

With basic properties of strongly rationally connected varieties
established, Theorem~\ref{theo:main} follows from the general result
(cf. \cite{kollar} IV.6.10.1):
\begin{theo}\label{theo:mainbis}
Let $\pi:\cY \ra B$ be a smooth morphism whose fibers are strongly rationally
connected.  Assume that $\pi$ has a section.
Then sections of $\cY\ra B$ satisfy approximation away from $S$.  
\end{theo}
\begin{proof}
Let $\overline{\pi}:\overline{\cY} \ra \overline{B}$ be a proper flat model of
$\cY\ra B$, which exists by \cite{Na}.  The section extends
to a section $\overline{s}$ of $\overline{\pi}$.  By a result of Artin and N\'eron
\cite{Art69b} Corollary 4.6, there exists a blow-up with center supported 
in $\overline{\pi}^{-1}(S)$
$$\widetilde{\cY} \ra \overline{\cY}$$ 
such that the proper transform of $\overline{s}(\overline{B})$ in $\widetilde{\cY}$
is contained in $\widetilde{\cY}^{sm}$.

Recall the proof of weak approximation at places of good reduction in Section 5 of \cite{HT}.
This is a bootstrap argument, using the existence of a section in the smooth locus
to construct sections with prescribed jets of successively higher order.  
{\em Properness} is used only to establish that the 
smooth fibers are strongly rationally connected, so we can produce very free
curves with desired properties.  In our situation, this 
is part of the hypotheses.
\end{proof}

Weak approximation at places of good reduction in positive characteristic was left
unresolved in \cite{HT}.  However, combining Theorem~\ref{theo:mainbis} with the 
main result of \cite{dJS} yields:
\begin{coro}\label{coro:positive}
Let $\pi:\cY \ra B$ be a smooth proper morphism with 
separably rationally connected fibers.  Then weak approximation
holds away from $S=\overline{B}\setminus B$.  
\end{coro}

\section{Cubic surfaces}
\label{sect:cubic}
We work over an algebraically closed field of characteristic zero.  
\begin{defi}
A {\em log Del Pezzo surface} is a pair $(X,\Delta)$ consisting
of a normal projective surface $X$ and an effective $\bQ$-divisor
$\Delta=\sum a_i \Delta_i, 0<a_i\le 1$ on $X$, with
log terminal singularities, such that $-(K_X+\Delta)$ is ample.
When $\Delta$ is empty, this is equivalent to saying that
$X$ has quotient singularities and ample anticanonical class.
\end{defi}

\begin{theo}[\cite{KeMc} 1.6]
\label{theo:KeMc}
The smooth locus of a log Del Pezzo surface $(X,\Delta)$ is
rationally connected, i.e., two generic points in $X^{sm}$
can be joined by an irreducible projective
rational curve contained in $X^{sm}$.  
\end{theo}
\begin{exam}[\cite{zhang}]
There exist projective rational surfaces with rational double points
whose smooth locus is not rationally connected.  Consider
$$\widetilde{X}=E\times \bP^1$$
where $(E,0)$ is an elliptic curve and the involution
\begin{eqnarray*}
\iota: \widetilde{X} & \ra & \widetilde{X} \\
        (e,[x_0,x_1]) & \mapsto & (-e,[x_1,x_0]). 
\end{eqnarray*}
The involution has eight isolated fixed points $\mfq \subset \widetilde{X}$.  
The quotient $X=\widetilde{X}/\left<\iota\right>$ has eight $A_1$ singularities
and is rational:  $X\ra E/\left<\iota\right>\simeq \bP^1$ is a conic bundle. 
Since $\widetilde{X}-\mfq \ra X^{sm}$ is a covering space, 
$\pi_1(X^{sm}) \subset \pi_1(\widetilde{X}-\mfq)$ with index two.  
Thus
$$\pi(\widetilde{X}-\mfq)\simeq \pi(\widetilde{X})\simeq \pi(E)\simeq \bZ\times \bZ$$
and $X^{sm}$ has infinite fundamental group.  However, rationally connected
varieties (even non-proper ones) have finite fundamental groups
(see Lemma 7.8 of \cite{KeMc} and Proposition 2.10 of \cite{KollarSM}, for example).  
\end{exam}

The following conjecture would allow us to 
apply Theorem~\ref{theo:main} to prove 
weak approximation for many log Del Pezzo surfaces:
\begin{conj} \label{conj:SRCDP}
The smooth locus of a log Del Pezzo surface is strongly rationally
connected.
\end{conj}

We prove this for cubic surfaces:
\begin{theo} \label{theo:cubicSRC}
Let $X\subset \bP^3$ be a cubic surface with rational double points.
Then $X^{sm}$ is strongly rationally connected.
\end{theo}
\begin{proof}
Let $x_1\in X^{sm}$ be a point.  We produce a rational curve
$R\subset X^{sm}$ joining $x_1$ and a generic point $x_2\in X^{sm}$.

We will make explicit precisely how $x_2$ must be chosen.  
We assume:
\begin{enumerate}
\item[(1)]{The tangent hyperplane section $H_2$ at $x_2$ 
is irreducible and nodal.}
\end{enumerate}
In particular, $H_2\subset X^{sm}$ and there are no lines 
$\ell\subset X$ containing $x_2$.
Projection from $x_2$ then gives a double cover
$$\mathrm{Bl}_{x_2}X \ra \bP^2;$$
the covering transformation interchanges 
the exceptional divisor and the proper transform.
We obtain a birational involution
\begin{eqnarray*}
\iota_{x_2}:X &\dashrightarrow& X \\
      x & \mapsto & x',
\end{eqnarray*}
where $\{x,x',x_2\}$ are collinear.
This factors as the blow-up of $x_2$ followed by
the blow-down of the proper transform of $H_2$.
Note that $\iota_{x_2}$ fixes
the singularities of $X$ and thus takes $X^{sm}$ to itself.

We also assume:
\begin{enumerate}
\item[(2)]{$H_2$
does not contain $x_1$.}
\end{enumerate}
It follows that 
$H_2$ does not contain $x'_1=\iota_{x_2}(x_1)$.
Moreover, $x_1$ and $x'_1$ are in the open subset on
which $\iota_{x_2}$ is an isomorphism. 

We assume furthermore:
\begin{enumerate}
\item[(3)]{$x_2$ is not contained in $H_1$.}
\end{enumerate}
It follows that $x_2 \not \in H'_1$,
the tangent hyperplane section at $x'_1$.  
Indeed, suppose that $x_2 \in H'_1$.  
We know that $x_2\neq x'_1$ (because $x'_1 \not \in H_2$),
so consider the line joining $x_2$ and $x'_1$.
This meets $X$ only
at $x_2$ and $x'_1$, so $x'_1=x_1$ and
$x_2 \in H_1$, a contradiction.  

Finally, we assume:
\begin{enumerate}
\item[(4)]{$H'_1$ is 
irreducible and nodal.}
\end{enumerate}
In particular, $H'_1 \subset X^{sm}$.  

Since $x_2 \not \in H'_1$, $\iota_{x_2}$ is regular along $H'_1$.  
We verify  that the rational curve
$R=\iota_{x_2}(H'_1)$ has the desired properties.
Since $H'_1\subset X^{sm}$ and
$\iota_{x_2}(X^{sm})\subset X^{sm}$, we find $R\subset X^{sm}$.   
We have $x'_1\in H'_1$, so $x_1=\iota_{x_2}(x'_1)\in R$.  
Since $H'_1$ meets $H_2$ in a point $y\neq x_2$, 
$x_2=\iota_{x_2}(y)\in R$.  
\end{proof}

We now prove Theorem~\ref{theo:cubic}:  For each singular fiber $\cX_b$,
$\cX_b^{sm}$ is strongly rationally connected by Theorem~\ref{theo:cubicSRC}.
Approximation follows from Theorem~\ref{theo:main}.

\begin{exam}
Here is another case where Conjecture~\ref{conj:SRCDP} is easily verified.
Let $X$ be a partial resolution of a cubic surface $\Sigma$ with at most $A_1$-singularities,
i.e., we have a factorization of the minimal resolution
$$\widetilde{\Sigma} \ra X \stackrel{\beta}{\ra} \Sigma.$$
Then $X^{sm}$ is strongly rationally connected.

Theorem~\ref{theo:cubic} implies that $\Sigma^{sm}$ is strongly rationally connected,
hence $\beta^{-1}(\Sigma^{sm})\subset (X^{sm})^0$.  The locus 
$X^{sm}\setminus \beta^{-1}(\Sigma^{sm})$ is a union of $(-2)$-curves $\{E_i\}$,
corresponding to the resolved singularities $\{p_i\}$ of $\Sigma$. 
If $(X^{sm})^0$ meets $E_i$, it must also contain $E_i$.  Hence
it suffices to show that for each $E_i$ there
exists a rational curve in
$X^{sm}$ meeting $E_i$ and  $\beta^{-1}(\Sigma^{sm})$ (see Remark~\ref{rema:include}).

To find this rational curve, consider the projection from $p_i$
$$\pi_i:\Sigma \dashrightarrow \bP^2$$
which induces a morphism $\pi'_i:X\ra \bP^2$.  The image of $E_i$ is a plane
conic and the image of the singularities of $X$
has codimension two in $\bP^2$, so there exists a rational curve 
$$f:\bP^1 \ra \bP^2\setminus \pi'_i(\mathrm{Sing}(X))$$
meeting the image of $E_i$.  

The same argument applies if $X$ is obtained from a cubic surface $\Sigma$ with
$A_1$ and $A_2$ singularities by resolving some subset of $\mathrm{Sing}(\Sigma)$.
\end{exam}

\bibliographystyle{plain}
\bibliography{sw}

\begin{thebibliography}{10}

\bibitem{Art69b}
M.~Artin.
\newblock Algebraic approximation of structures over complete local rings.
\newblock {\em Inst. Hautes \'Etudes Sci. Publ. Math.}, (36):23--58, 1969.

\bibitem{Br1}
E.~Brieskorn.
\newblock Die {A}ufl\"osung der rationalen {S}ingularit\"aten holomorpher
  {A}bbildungen.
\newblock {\em Math. Ann.}, 178:255--270, 1968.

\bibitem{Br2}
E.~Brieskorn.
\newblock Singular elements of semi-simple algebraic groups.
\newblock In {\em Actes du Congr\`es International des Math\'ematiciens (Nice,
  1970), Tome 2}, pages 279--284. Gauthier-Villars, Paris, 1971.

\bibitem{CT}
J.~L. Colliot-Th{\'e}l{\`e}ne and P.~Gille.
\newblock Remarques sur l'approximation faible sur un corps de fonctions d'une
  variable.
\newblock In {\em Arithmetic of higher-dimensional algebraic varieties (Palo
  Alto, CA, 2002)}, volume 226 of {\em Progr. Math.}, pages 121--134.
  Birkh\"auser Boston, Boston, MA, 2004.

\bibitem{dJS}
A.~J. de~Jong and J.~Starr.
\newblock Every rationally connected variety over the function field of a curve
  has a rational point.
\newblock {\em Amer. J. Math.}, 125(3):567--580, 2003.

\bibitem{De}
O.~Debarre.
\newblock Vari\'et\'es rationnellement connexes (d'apr\`es {T}. {G}raber, {J}.
  {H}arris, {J}. {S}tarr et {A}. {J}. de {J}ong).
\newblock {\em Ast\'erisque}, 290:Exp. No. 905, ix, 243--266, 2003.
\newblock S\'eminaire Bourbaki. Vol. 2001/2002.

\bibitem{GHS}
T.~Graber, J.~Harris, and J.~Starr.
\newblock Families of rationally connected varieties.
\newblock {\em J. Amer. Math. Soc.}, 16(1):57--67, 2003.

\bibitem{HT}
B.~Hassett and Y.~Tschinkel.
\newblock Weak approximation over function fields.
\newblock {\em Inventiones mathematicae}, 163:171--190, 2006.

\bibitem{KeMc}
S.~Keel and J.~McKernan.
\newblock Rational curves on quasi-projective surfaces.
\newblock {\em Mem. Amer. Math. Soc.}, 140(669):viii+153, 1999.

\bibitem{KollarSM}
J.~Koll{\'a}r.
\newblock {\em Shafarevich maps and automorphic forms}.
\newblock M. B. Porter Lectures. Princeton University Press, Princeton, NJ,
  1995.

\bibitem{kollar}
J.~Koll{\'a}r.
\newblock {\em Rational curves on algebraic varieties}, volume~32 of {\em
  Ergebnisse der Math.}
\newblock Springer-Verlag, Berlin, 1996.

\bibitem{KMM}
J.~Koll{\'a}r, Y.~Miyaoka, and S.~Mori.
\newblock Rationally connected varieties.
\newblock {\em J. Algebraic Geom.}, 1(3):429--448, 1992.

\bibitem{Ma}
D.~Madore.
\newblock Approximation faible aux places de bonne r\'eduction sur les surfaces
  cubiques sur les corps des fonctions, 2004.
\newblock preprint.

\bibitem{Na}
M.~Nagata.
\newblock A generalization of the imbedding problem of an abstract variety in a
  complete variety.
\newblock {\em J. Math. Kyoto Univ.}, 3:89--102, 1963.

\bibitem{SWD}
P.~Swinnerton-Dyer.
\newblock Weak approximation and {$R$}-equivalence on cubic surfaces.
\newblock In {\em Rational points on algebraic varieties}, volume 199 of {\em
  Progr. Math.}, pages 357--404. Birkh\"auser, Basel, 2001.

\bibitem{zhang}
D.-Q. Zhang.
\newblock Algebraic surfaces with nef and big anti-canonical divisor.
\newblock {\em Math. Proc. Cambridge Philos. Soc.}, 117(1):161--163, 1995.

\end{thebibliography}

\end{document}